\documentstyle[12pt,amsfonts,leqno,amscd]{amsart}

\newtheorem{Theorem}{Theorem}[section]

\newtheorem{Proposition}[Theorem]{Proposition}

\newtheorem{Lemma}[Theorem]{Lemma}

\theoremstyle{remark}
\newtheorem{Example}[Theorem]{Example}

\begin{document}
\title{Entire pluricomplex Green functions and Lelong numbers of projective currents}
\author{Dan Coman}
\thanks{The author was supported by the NSF grant DMS 0140627}
\address{dcoman@@syr.edu, Department of Mathematics,  215 Carnegie Hall, Syracuse
University, Syracuse, NY 13244-1150, USA}

\begin{abstract}Let $T$ be a positive closed current of bidimension (1,1) and
unit mass on the complex projective space ${\Bbb P}^n$. We prove
that the set $V_\alpha(T)$ of points where $T$ has Lelong number
larger than $\alpha$ is contained in a complex line if
$\alpha\geq2/3$, and $|V_\alpha(T)\setminus L|\leq1$ for some
complex line $L$ if $1/2\leq\alpha<2/3$. We also prove that in
dimension 2 and if $2/5\leq\alpha<1/2$, then
$|V_\alpha(T)\setminus C|\leq1$ for some conic $C$.
\end{abstract}
\maketitle

\section{Introduction}
\par Let $T$ be a positive closed current of bidimension (1,1) on
the complex projective space ${\Bbb P}^n$, which has mass
$\|T\|=1$ (see Section \ref{S:green} for the defintion of
$\|T\|$). Siu's theorem \cite{Siu74} states that the upper level
sets $E_\alpha(T)$ of the Lelong numbers $\nu(T,\cdot)$ of $T$,
$$E_\alpha(T):=\{p\in{\Bbb P}^n:\,\nu(T,p)\geq\alpha\},$$ are
analytic subvarieties of ${\Bbb P}^n$ of dimension at most one.
Hence by Chow's theorem, $E_\alpha(T)$ are algebraic varieties.
\par We consider here the following upper level sets of Lelong
numbers of $T$:
$$V_\alpha(T):=\{p\in{\Bbb P}^n:\,\nu(T,p)>\alpha\},\;\alpha<1.$$
Given a finite set $S$, we denote by $|S|$ the number of points of
$S$. The main results of this paper are the following theorems.
\begin{Theorem}\label{T:234} (i) If $\alpha\geq2/3$, the set
$V_\alpha(T)$ is either a complex line or an at most countable
subset of a complex line.
\par (ii) If $1/2\leq\alpha<2/3$, the set $V_\alpha(T)$ is either
a complex line, or a countable subset of a complex line, or a
finite set such that $|V_\alpha(T)\setminus L|\leq1$ for some
complex line $L$.\end{Theorem}
\begin{Theorem}\label{T:25} If $n=2$ and $2/5\leq\alpha<1/2$, the set
$V_\alpha(T)$ is one of the following: a conic, the union of a
complex line and an at most countable subset of a complex line, a
countable subset of a conic, a finite set such that
$|V_\alpha(T)\setminus C|\leq1$ for some conic $C$.
\end{Theorem}
\par The proofs of these theorems require the construction of
entire plurisubharmonic functions of logarithmic growth and with
logarithmic poles in some special finite subsets of ${\Bbb C}^n$.
In certain cases such functions were constructed in \cite{Co} and
\cite{CN}.
\par If $D\subseteq{\Bbb C}^n$ is an open set, we denote by $PSH(D)$
the class of plurisubharmonic functions in $D$. A function $u\in
PSH(D)$ is said to have a logarithmic pole of weight $\alpha>0$ at
$p\in D$, if
$$u(z)=\alpha\log\|z-p\|+O(1)$$ holds for $z\in D$ near $p$. We say that
$u\in PSH({\Bbb C}^n)$ has logarithmic growth if
$$\gamma_u:=\limsup_{\|z\|\rightarrow+\infty}\frac{u(z)}{\log\|z\|}<+\infty.$$
\par Let $S$ be a finite subset of ${\Bbb C}^n$. If $u\in PSH({\Bbb C}^n)$
has logarithmic growth, it is locally bounded and maximal on
${\Bbb C}^n\setminus S$, and it equals $-\infty$ on $S$, we call
$u$ an entire pluricomplex Green function of $S$. We let $E(S)$ be
the class of entire pluricomplex Green functions of $S$ with
logarithmic poles of weight 1 at all points of $S$. Moreover, we
let $\widetilde E(S)\subset PSH({\Bbb C}^n)$ be the class of
plurisubharmonic functions with logarithmic growth, which are
locally bounded on ${\Bbb C}^n\setminus S$ and have logarithmic
poles of weight 1 at the points of $S$. These classes were
introduced in \cite{CN}, where we defined and studied two affine
invariants $\gamma(S)\geq\widetilde\gamma(S)$ of the set $S$,
$$\gamma(S)=\inf\{\gamma_u:\,u\in E(S)\},\;
\widetilde\gamma(S)=\inf\{\gamma_u:\,u\in\widetilde E(S)\}.$$
These numbers are connected to the singular degree of $S$
introduced in \cite{W77} and studied in \cite{W77} and \cite{Ch}.
\par We note that entire pluricomplex Green functions were
generally considered in the case when they are locally bounded,
for example when dealing with the extremal function of a compact
set, while pluricomplex Green functions with logarithmic poles
were considered mostly in the case of bounded domains, by
prescribing their boundary values to be 0 (see e.g. the survey
\cite{B}). For our purposes, it is useful to combine these two
cases, and consider entire plurisubharmonic functions with
finitely many logarithmic poles, as in \cite{Co} and \cite{CN}.
\par Theorems \ref{T:234} and \ref{T:25}, together with two other
related results, are proved in Section \ref{S:lel}. We also give
examples of currents which show that the results of these theorems
are sharp. In Proposition \ref{P:lel} from Section \ref{S:green}
we show the connection between the Lelong numbers of projective
currents at the points of a finite set $S$, and the growth
constants $\gamma_u$ of functions $u\in PSH({\Bbb C}^n)$ with
logarithmic poles in $S$. The entire pluricomplex Green functions
which we need to prove the results from Section \ref{S:lel} are
constructed in Section \ref{S:green}.

\section{Entire pluricomplex Green functions}\label{S:green}
\par We denote by $[z:t]$, $(z,t)\in{\Bbb C}^n\times{\Bbb C}$,
$(z,t)\neq(0,0)$, the homogeneous coordinates on ${\Bbb P}^n$, and
we use the standard embedding ${\Bbb C}^n\equiv\{[z:1]\in{\Bbb
P}^n:\,z\in{\Bbb C}^n\}$. Let $\omega$ be the standard K\"ahler
form on ${\Bbb P}^n$, corresponding to the Fubini-Study metric.
Then
$$\omega\mid_{{\Bbb C}^n}=dd^cV,\;V(z):=\log\sqrt{1+\|z\|^2},$$
where $d^c=(\partial-\overline\partial)/(2\pi i)$. If $T$ is a
positive closed current on ${\Bbb P}^n$ of bidimension $(l,l)$,
its mass is given by
$$\|T\|=\int_{{\Bbb P}^n}T\wedge\omega^l.$$
\par The following simple result is analogous to \cite[Theorem
3.4]{CN}.
\begin{Proposition}\label{P:lel}Let $S=\{p_1,\dots,p_k\}\subset{\Bbb C}^n$
and let $T$ be a positive closed current on ${\Bbb P}^n$ of
bidimension $(l,l)$. If $u\in PSH({\Bbb C}^n)$ has logarithmic
growth, it is locally bounded outside a finite set, and
$u(z)\leq\alpha_j\log\|z-p_j\|+O(1)$ for $z$ near $p_j$, where
$\alpha_j>0$, $1\leq j\leq k$, then
$$\sum_{j=1}^k\alpha_j^l\nu(T,p_j)\leq\gamma_u^l\|T\|.$$ In
particular,
$\sum_{j=1}^k\nu(T,p_j)\leq\widetilde\gamma(S)^l\|T\|$.
\end{Proposition}
\begin{pf} Since $u$ is locally bounded outside a finite set, the
positive measure $T\wedge(dd^cu)^l$ is well defined on ${\Bbb
C}^n$, by \cite[Proposition 2.1]{D93}. Demailly's first comparison
theorem for Lelong numbers with weights \cite[Theorem 5.1]{D93}
implies that
$$T\wedge(dd^cu)^l(\{p_j\})\geq\alpha_j^l\nu(T,p_j).$$ Using
\cite[Proposition 3.2]{CN} we obtain
$$\sum_{j=1}^k\alpha_j^l\nu(T,p_j)\leq\int_{{\Bbb C}^n}T\wedge(dd^cu)^l\leq
\gamma_u^l\int_{{\Bbb C}^n}T\wedge(dd^cV)^l\leq\gamma_u^l\|T\|.$$\end{pf}
\par We now construct entire pluricomplex Green functions of small
growth, for special finite subsets $S$ of ${\Bbb C}^n$. We denote
by $m_j(S)$ the maximum number of points of $S$ which are
contained in an algebraic curve of degree $j$.
\begin{Lemma}\label{L:line} If $S\subset{\Bbb C}^n$, $|S|\in\{3,4\}$,
$m_1(S)=2$, then $\gamma(S)=2$.\end{Lemma}
\begin{pf} If $S$ is contained in a two-dimensional complex plane,
the lemma follows from analogous results in ${\Bbb C}^2$
\cite{CN}. Otherwise, we can assume that $S\subset{\Bbb
C}^3\times\{0\}$ consists of the points $(0,0,0,0),\,(1,0,0,0),\,
(0,1,0,0),\,(0,0,1,0)$. Let
$$u(z)=\frac{1}{2}\,
\log\left(\sum_{j=1}^3|P_j(z_1,z_2,z_3)|^2+\sum_{k=4}^n|z_k|^2\right),
\;z=(z_1,\dots,z_n),$$ where
$P_j(z_1,z_2,z_3)=z_j(z_1+2z_2+3z_3-j)$. Then
$\gamma_u=\gamma(S)=2$.\end{pf}
\par In the next three results, we consider the case of sets $S\subset{\Bbb
C}^2$ with 7 or 8 elements.
\begin{Proposition}\label{P:conic1} Let $S\subset{\Bbb C}^2$ be such
that $|S|=7$ and $m_2(S)=5$. Then $S$ has an entire pluricomplex
Green function $u$ with $\gamma_u=4$, such that $u$ has
logarithmic poles of weight 2 at 3 of the points of $S$, and of
weight 1 at the remaining 4 points of $S$.\end{Proposition}
\begin{pf} We show that we can find 3 points of $S$, $\zeta_1,\zeta_2,\zeta_3$,
with the following property: There exist two polynomials $P_1,P_2$
of degree 4, with no common factors, such that $P_1(S)=P_2(S)=0$
and both $P_1,P_2$ vanish to second order at each point $\zeta_j$.
By Bezout's theorem it follows that $S=\{P_1=P_2=0\}$, and the
intersection numbers $(P_1\cdot P_2)_{\zeta_j}=4$ and $(P_1\cdot
P_2)_\zeta=1$ at all other points $\zeta\in S$. By \cite[Theorem
4.1]{CN}, the function $u=\frac{1}{2}\,\log(|P_1|^2+|P_2|^2)$ has
the desired properties.
\par Let $S=\{p_1,\dots,p_7\}$. The vector space of polynomials of
degree at most 4, which vanish on $S$ and to second order at 3
given points of $S$, has dimension at least 2.
\par {\em Case 1.} $m_1(S)=2,\;m_2(S)=5$. Let $C_1,C_2$ be quadratic
polynomials vanishing at $p_1,p_2,p_3,p_4,p_5$, respectively at
$p_1,p_2,p_3,p_6,p_7$. Then $C_j$ are irreducible, the conics
$C_j=0$ are smooth algebraic curves, $C_1$ does not vanish at
$p_6,p_7$, and $C_2$ does not vanish at $p_4,p_5$. Let
$P_1=C_1C_2$ and $P_2$ be a polynomial of degree 4 vanishing on
$S$ and to second order at $p_1,p_2,p_3$, such that $P_1,P_2$ are
linearly independent. If $C_1$ divides $P_2$, then $P_2=C_1C$,
where $C$ is a quadratic polynomial vanishing at
$p_1,p_2,p_3,p_6,p_7$, so $C=\alpha C_2$. It follows that
$P_1,P_2$ have no common factors.
\par {\em Case 2.} $m_1(S)=3,\;m_2(S)=5$. We may assume that
$p_1,p_2,p_3$ lie on a complex line $l$, hence $p_4,p_5,p_6,p_7$
lie outside of $l$ and are in general position, since $m_2(S)=5$.
\par (i) Assume first that $\{p_1,p_2,p_3\}$ is the only 3 point
subset of $S$ contained in a complex line. Let $C_1,C_2$ be
irreducible quadratic polynomials vanishing at
$p_4,p_5,p_6,p_1,p_2$, and respectively at $p_4,p_5,p_6,p_3,p_7$.
We let $P_1=C_1C_2$ and continue as in Case 1.
\par (ii) Assume that, after relabelling if necessary the points of
$S$, $p_1,p_4,p_5$ are also contained in a complex line, so
$p_2,p_3,p_6,p_7$ are in general position. We show that there
exist 3 points of $S$, $p_i,p_j,p_k$, with the following property
$(P)$: $p_i,p_j,p_k$ are not contained in a complex line and the
polynomial $L_{ij}L_{jk}L_{ki}$ does not vanish at any of the
remaining 4 points of $S$. Then we can continue as in (i), with
the points $p_4,p_5,p_6$ replaced by $p_i,p_j,p_k$ (although the
corresponding conics may now be reducible).
\par If $(L_{24}L_{25}L_{34}L_{35})(p_j)\neq0$ for $j=6,7$, then
the points $p_2,p_4,p_6$ have the above property. Otherwise, we
may assume that $L_{24}(p_6)=0$. Since $m_2(S)=5$, we have
$L_{35}(p_7)\neq0$. If $L_{35}(p_6)\neq0$ then $p_3,p_5,p_6$
verify property $(P)$. Finally, we assume
$L_{24}(p_6)=L_{35}(p_6)=0$. If $L_{25}(p_7)\neq0$ or
$L_{34}(p_7)\neq0$, then the points $p_2,p_5,p_7$, respectively
$p_3,p_4,p_7$, have property $(P)$. If
$L_{25}(p_7)=L_{34}(p_7)=0$, then $p_1,p_6,p_7$ verify property
$(P)$. \end{pf}
\begin{Proposition}\label{P:conic2} Let $A\subset{\Bbb C}^2$ with
$|A|=7$, $m_1(A)\leq3$, $m_2(A)=6$, and let $\Gamma$ be the conic
with $|A\cap\Gamma|=6$. Let $q\not\in A\cup\Gamma$. Then
$m_1(A\cup\{q\})\leq4$ and we have:
\par (i) If $m_1(A\cup\{q\})\leq3$ there exists $u\in PSH({\Bbb C}^2)$
with $\gamma_u=3$, such that $u$ is locally bounded outside a
finite set and $u(z)\leq\log\|z-p\|+O(1)$ near each point $p\in
A\cup\{q\}$.
\par (ii) If $m_1(A\cup\{q\})=4$ there exists a subset $S$
of $A\cup\{q\}$ with 7 elements, which has an entire pluricomplex
Green function as in the conclusion of Proposition
\ref{P:conic1}.\end{Proposition}
\begin{pf} Let $\{p_1\}=A\setminus\Gamma$,
$\{p_2,\dots,p_7\}=A\cap\Gamma=A'$, and $C$ be a quadratic
polynomial defining $\Gamma$. We have either $m_1(A')=2$ and $C$
is irreducible, or else $C=l_1l_2$, where $l_j$ are linear
polynomials and $|\{l_j=0\}\cap A'|=3$,
$\{l_1=0\}\cap\{l_2=0\}\cap A=\emptyset$. In the latter case, we
can assume that $p_2,p_4,p_6$ lie on the line $l_1=0$ and
$p_3,p_5,p_7$ on the line $l_2=0$. In either case, the conic
$\Gamma$ is smooth at the points $p_j$, $j\geq2$. \vspace{2mm}\par
$(i)$ If $m_1(A\cup\{q\})\leq3$, let $L$ be a linear polynomial
with $L(p_1)=L(q)=0$, and let $P_1=LC$. There exists a polynomial
$P_2$ of degree 3 which is zero on $A\cup\{q\}$ and such that
$P_1,P_2$ are linearly independent. Assume that $P_2=LC'$, for
some quadratic polynomial $C'$. Since $L$ vanishes at most at one
of the points $p_2,\dots,p_7$, and the remaining 5 points
determine $\Gamma$ uniquely, it follows that $C'=\alpha C$, a
contradiction. Similarly one shows that $C$ (or $l_j$ if $C$ is
reducible) cannot divide $P_2$, so $P_1,P_2$ have no common
factors.  We let $u=\frac{1}{2}\,\log(|P_1|^2+|P_2|^2)$.
\vspace{2mm}\par $(ii)$ If $m_1(A\cup\{q\})=4$ then the complex
line determined by $p_1,q$ must contain two points of $\Gamma$.
After relabelling points, we can assume that $p_1,p_2,p_3,q$ are
contained in a complex line. Let $L_{jk}$ be a linear polynomial
such that $L_{jk}(p_j)=L_{jk}(p_k)=0$. Since $m_1(A)\leq3$, we
have for each $j\geq2$ that $L_{1j}$ can vanish at most at one
other point $p_k$ ($k\geq2$, $k\neq j$).
\par{\em Case 1.} If some polynomial $L_{1j}$, $j\geq4$,
does not vanish at any other point $p_k$, then there must exist
another polynomial $L_{1i}$, $i\geq4$, with the same property.
Note that when $C=l_1l_2$ then such $i$ exists so that $p_i,p_j$
are not both contained in a line $l_k=0$. We let $S=A$,
$P_1=L_{1i}L_{1j}C$, and $P_2$ be a polynomial of degree 4
vanishing on $S$ and to second order at $p_1,p_i,p_j$, such that
$P_1,P_2$ are linearly independent. A direct application of
Bezout's theorem shows that $P_1,P_2$ have no common factors. By
the considerations at the beginning of the proof of Proposition
\ref{P:conic1}, $u=\frac{1}{2}\,\log(|P_1|^2+|P_2|^2)$ is the
desired pluricomplex Green function of $S$.
\par{\em Case 2.} If each polynomial $L_{1j}$, $j\geq4$, vanishes
at some other point $p_k$, $k\geq4$, then after relabelling points
we have $L_{14}(p_5)=L_{16}(p_7)=0$. If
$S=(A\cup\{q\})\setminus\{p_2\}$ then $m_2(S)=5$, and the
conclusion follows from Proposition \ref{P:conic1}. Indeed, the
points $p_3,\dots,p_7$ determine $\Gamma$, which does not contain
$p_1,q$, and $p_1,p_4,p_5,p_6,p_7$ determine the conic
$L_{45}L_{67}=0$, which does not contain $p_3,q$. If
$S\setminus\{p_j\}$, $j\geq4$, is contained in a conic and without
loss of generality $j=4$, then this conic is defined by
$L_{23}L_{67}=0$. But $L_{23}L_{67}(p_5)\neq0$.\end{pf}

\section{Upper level sets of Lelong numbers}\label{S:lel}
\par Throughout this section, $T$ is a positive closed current
of bidimension $(1,1)$ on ${\Bbb P}^n$, normalized by $\|T\|=1$.
Then the Lelong numbers $\nu(T,p)\leq1$, at every point $p\in{\Bbb
P}^n$. We will need the following simple lemma:
\begin{Lemma}\label{L:count} If $S$ is an at most countable subset
of ${\Bbb P}^n$, then there exists a hyperplane $H$ which does not
intersect $S$.\end{Lemma}
\begin{pf} We show by induction on $k$, $0\leq k\leq n-1$, that
there exists a $k$-dimensional plane $P_k$ which does not
intersect $S$. This is clear for $k=0$, and assume such $P_k$
exists for $k<n-1$. The family of $(k+1)$-dimensional planes $P$
which contain $P_k$ is uncountable, and the sets $P\setminus P_k$
are pairwise disjoint. Since $S$ is at most countable and it does
not intersect $P_k$, it follows that there is a
$(k+1)$-dimensional plane $P_{k+1}\supset P_k$ which does not
intersect $S$.\end{pf}
\par Recall that for a finite subset $S$ of ${\Bbb P}^n$, $m_j(S)$
denotes the maximum number of points in $S$ which are contained in
an algebraic curve of degree $j$. We now proceed with the proofs
of Theorems \ref{T:234} and \ref{T:25}.
\par\vspace{3mm}\noindent{\em Proof of Theorem \ref{T:234}.}
$(i)$ Assume that some set $E_\beta(T)$, $\beta>\alpha$, has
dimension one, so it contains an algebraic curve $C$. Then
$T=\beta[C]+R$, where $[C]$ is the current of integration along
$C$ and $R$ is a $(1,1)$ bidimensional positive closed current on
${\Bbb P}^n$. The degree of $C$ (see e.g. \cite{LG}) is
$\|[C]\|\leq1/\beta<2$, so $C$ is a complex line. Since
$\nu(R,p)\leq\|R\|=1-\beta<1/3$ at all $p\in{\Bbb P}^n$, it
follows that $V_\alpha(T)=C$.
\par If all sets $E_\beta(T)$, $\beta>\alpha$, have dimension 0,
then $V_\alpha(T)$ is at most countable. By Lemma \ref{L:count}
there is a hyperplane $H$ so that $V_\alpha(T)\subset{\Bbb
P}^n\setminus H$, so we may assume $V_\alpha(T)\subset{\Bbb C}^n$.
If $V_\alpha(T)$ is not contained in a complex line, there exists
$S\subseteq V_\alpha(T)$ with $|S|=3$ and $m_1(S)=2$, so
$\gamma(S)=2$ by Lemma \ref{L:line}. Proposition \ref{P:lel}
yields the following contradiction: $$3\alpha<\sum_{p\in
S}\nu(T,p)\leq\gamma(S)\|T\|=2.$$
\par $(ii)$ Arguing as in the proof of $(i)$, we have either
that $V_\alpha(T)$ is a complex line, or it is at most countable
and contained in ${\Bbb C}^n$.
\par If $V_\alpha(T)$ is countable, then $E_{1/2}(T)$ must contain
an algebraic curve $C$. As before, $T=\frac{1}{2}\,[C]+R$, so the
degree of $C$ is at most 2. If $C$ has degree 2 then $R=0$, so
$\nu(T,p)>1/2$ only at the singular points of $C$ and $V_{1/2}(T)$
is a finite set. We conclude that $C$ is a complex line. Since
$\|R\|=1/2$, we have $V_\alpha(T)\subset C$.
\par We assume finally that $V_\alpha(T)$ is a finite set not
contained in a complex line and with at least 4 elements. So there
exists $S=\{p_1,p_2,p_3,p_4\}\subset V_\alpha(T)$ such that
$p_1,p_2,p_3$ are not contained in a complex line. If $m_1(S)=2$
then Lemma \ref{L:line} and Proposition \ref{P:lel} imply
$4\alpha<2$, a contradiction. So, after relabelling points, $p_4$
lies on the complex line $L$ determined by $p_1,p_2$. If there
exists $p\in V_\alpha(T)\setminus L$, $p\neq p_3$, then at least
one of the following sets $S$,
$$\{p_1,p_2,p_3,p\},\;\{p_1,p_4,p_3,p\},\;\{p_2,p_4,p_3,p\},$$ has
$m_1(S)=2$, and we reach a contradiction as above. Therefore
$|V_\alpha(T)\setminus L|=1$. $\Box$ \par\vspace{3mm}\noindent{\em
Proof of Theorem \ref{T:25}.} As in the previous proof, there are
two possibilities: some set $E_\beta(T)$, $\beta>\alpha$, contains
an algebraic curve $C$, or $V_\alpha(T)$ is at most countable and
contained in ${\Bbb C}^2$.
\par In the first case, $T=\beta[C]+R$, for some positive closed
(1,1) current $R$ on ${\Bbb P}^2$, and the degree of $C$ is at
most 2. If $C$ is a conic then $\|R\|<1/5$ and $V_\alpha(T)=C$. If
$C$ is a complex line and $R=0$ then $V_\alpha(T)=C$. Otherwise
$0<\|R\|=1-\beta<3/5$ and $V_\alpha(T)=C\cup V_\alpha(R)$. Since
$\alpha/\|R\|\geq2/3$, $V_\alpha(R)$ is a complex line or an at
most countable subset of a complex line, by Theorem \ref{T:234}.
\par We assume next that $V_\alpha(T)$ is countable, so $E_{2/5}(T)$
contains an algebraic curve $C$ of degree at most 2, and
$T=\frac{2}{5}\,[C]+R$. If $C$ is a conic then $V_\alpha(T)\subset
C$. If $C$ is a complex line then $\|R\|=3/5$, so $V_\alpha(R)$ is
contained in a complex line since $\alpha/\|R\|\geq2/3$, and
$V_\alpha(T)\subset C\cup V_\alpha(R)$.
\par We assume finally that $V_\alpha(T)$ is a finite set and
$|V_\alpha(T)\setminus C|>1$, for any conic $C$. Then
$V_\alpha(T)$ has a subset $A$ with $|A|=7$. We have the following
possibilities:
\par {\em Case 1.} $m_2(A)=5$. Then $A$ has an entire pluricomplex
Green function $u$ as in Proposition \ref{P:conic1}. Proposition
\ref{P:lel} applied to $u$ and $T$ implies $10\alpha<4$, a
contradiction.
\par {\em Case 2.} $m_1(A)\leq3,\;m_2(A)=6$. Let $C$ be the conic
with $|A\cap C|=6$. Since $|V_\alpha(T)\setminus C|>1$, there is
$q\in V_\alpha(T)\setminus(A\cup C)$. Applying Proposition \ref
{P:lel} with the functions provided by Proposition \ref{P:conic2}
yields $8\alpha<3$, or $10\alpha<4$, a contradiction.
\par {\em Case 3.} $m_1(A)=2$ and $A$ is contained in a conic $C$.
There exists $q\in V_\alpha(T)\setminus C$. If $p\in A$, then
$S=(A\cup\{q\})\setminus\{p\}$ verifies the hypotheses of Case 2.
\par {\em Case 4.} $m_1(A)\geq3$. We show that $V_\alpha(T)$ has
a subset $S$ with $|S|=7$, $m_1(S)\leq3$, $m_2(S)\leq6$, so we are
back to Case 1 or Case 2. Let $L_{jk}$ be the complex line
determined by the points $p_j,p_k$. Let $p_1,p_2,p_3$ be points of
$A$ so that $p_3\in L_{12}$. As $|V_\alpha(T)\setminus C|>1$ for
any conic $C$, there exist points $p_4,\dots,p_7\in
V_\alpha(T)\setminus L_{12}$, with $p_6,p_7\not\in L_{45}$. Let
$S_1=\{p_1,\dots,p_7\}$. Then $m_2(S_1)\leq6$, so $m_1(S_1)\leq4$.
If $m_1(S_1)\leq3$ we let $S=S_1$. If $m_1(S_1)=4$ then, after
relabelling points, $p_6,p_7\in L_{14}$. Hence there exists
$p_8\in V_\alpha(T)\setminus(L_{12}\cup L_{14})$, $p_8\neq p_5$.
The point $p_8$ can lie on at most one of the lines $L_{5j}$,
$j=4,6,7$. If $p_8\in L_{5j}$ for some $j\in\{4,6,7\}$, we let
$S=(S_1\cup\{p_8\})\setminus\{p_j\}$. Otherwise, we let
$S=(S_1\cup\{p_8\})\setminus\{p_7\}$. $\Box$
\par It is not difficult to see that all cases described in Theorems
\ref{T:234} and \ref{T:25} can occur. We present here a few
examples in ${\Bbb P}^2$.
\begin{Example}\label{E:23s} Let $L_1,L_2,L_3$ be complex lines with
$L_1\cap L_2\cap L_3=\emptyset$, and let $\{p_1\}=L_2\cap L_3$,
$\{p_2\}=L_1\cap L_3$, $\{p_3\}=L_1\cap L_2$. If
$T=\frac{1}{3}\,\sum_{j=1}^3[L_j]$ then
$E_{2/3}(T)=\{p_1,p_2,p_3\}$ is not contained in a complex line,
so the value $\alpha=2/3$ in Theorem \ref{T:234} $(i)$ is the best
possible.
\end{Example}
\begin{Example}\label{E:12s} Let $L_1,\dots,L_4$ be complex lines
so that no 3 of them pass through the same point, let
$\{p_{jk}\}=L_j\cap L_k$, $1\leq j<k\leq4$, and $S=\{p_{jk}\}$,
$|S|=6$. If $T=\frac{1}{4}\,\sum_{j=1}^4[L_j]$ then $E_{1/2}(T)=S$
has $m_1(S)=3$, so the value $\alpha=1/2$ in Theorem \ref{T:234}
$(ii)$ is sharp. Moreover, $E_{1/2}(T)$ is not contained in a
conic.
\end{Example}
\begin{Example}\label{E:12c} If $C$ is a conic and $T=[C]/2$ then
$E_{1/2}(T)\setminus L$ is an uncountable set, for every complex
line $L$.\end{Example}
\begin{Example}\label{E:12f} Let $L_j$ be complex lines containing
the point $q$, $L$ be a complex line with $q\not\in L$, and
$\{p_j\}=L\cap L_j$. If $m\geq2$ and
$$T=\frac{m-1}{2m}\,[L]+\frac{m+1}{2m^2}\,\sum_{j=1}^m[L_j],$$
then $V_{1/2}(T)=\{p_1,\dots,p_m,q\}$ has $m$ points on $L$.
\end{Example}
\begin{Example}\label{E:25s} Let $L_j,p_j$ be as in Example \ref{E:23s},
and let $p_4\not\in L_1\cup L_2\cup L_3$. Let $L_4,L_5,L_6$ be the
complex lines determined by the points $p_4$ and respectively
$p_1,p_2,p_3$, and let $\{p_5\}=L_1\cap L_4$, $\{p_6\}=L_2\cap
L_5$, $\{p_7\}=L_3\cap L_6$. Finally, let $l_1,l_2,l_3$ be the
complex lines determined by pairs of points from
$\{p_5,p_6,p_7\}$. If $S=\{p_1,\dots,p_7\}$ and
$$T=\frac{1}{15}\left(2\sum_{j=1}^6[L_j]+\sum_{j=1}^3[l_j]\right),$$
then $m_2(S)=5$, $E_{2/5}(T)=S$, so the value $\alpha=2/5$ in
Theorem \ref{T:25} is sharp.
\end{Example}
\begin{Example}\label{E:25count} Let $C$ be a conic and let $p_j\in C$,
$j\geq0$, be distinct points so that $p_j\rightarrow p_0$ as
$j\rightarrow\infty$. Let $L_j$ be a complex line passing through
$p_j$ and some given point $q\not\in C$. If $\epsilon_j>0$,
$\sum_{j=0}^\infty\epsilon_j=1/5$, and
$$T=\frac{2}{5}\,[C]+\sum_{j=0}^\infty\epsilon_j[L_j],$$ then
$V_{2/5}(T)$ is a countable subset of $C$.
\end{Example}
\par We saw in Examples \ref{E:12s} and \ref{E:12c} that for
$\beta<1/2$ one can have $|V_\beta(T)\setminus L|>1$ for every
complex line $L$. But if $T$ has ``large" Lelong numbers at two
points, then the conclusion of Theorem \ref{T:234} $(ii)$ still
holds for some values $\beta<1/2$.
\begin{Theorem}\label{T:magnet} Assume that $\alpha>1/2$ and there
are points $q_1,q_2\in{\Bbb P}^n$ so that $\nu(T,q_j)\geq\alpha$,
$j=1,2$. If $\beta=(2-\alpha)/3$ then $|V_\beta(T)\setminus
L|\leq1$ for some complex line $L$ which contains at least one of
the points $q_1,q_2$.\end{Theorem}
\begin{pf} Let $L_1$ be the complex line determined by $q­_1,q_2$.
If $|V_\beta(T)\setminus L_1|>1$, let $p_1,p_2\in
V_\beta(T)\setminus L_1$ and let $L_2$ be the complex line
determined by $p_1,p_2$. We choose $1/2<\alpha'<\alpha$ so that
$\nu(T,p_j)>(2-\alpha')/3$ and we consider the current
$$R=\frac{2\alpha'-1}{2\alpha'}\,[L_2]+\frac{1}{2\alpha'}\,T.$$ Then
$\nu(R,q_j)\geq\alpha/(2\alpha')>1/2$ and
$$\nu(R,p_j)>\frac{2\alpha'-1}{2\alpha'}+\frac{2-\alpha'}{6\alpha'}
>\frac{1}{2}\;,$$ for $j=1,2$. Theorem \ref{T:234} $(ii)$ implies
that one of the points $q_1,q_2$, say without loss of generality
$q_1$, lies on $L_2$. We show that
$V_\beta(T)\setminus\{q_2\}\subseteq L_2$. If not, there exists
$p_3\in V_\beta(T)\setminus\{q_2\}$, $p_3\not\in L_2$. If $L_{jk}$
denotes the complex line determined by $p_j,p_k$, let
$$R=\frac{2\alpha'-1}{6\alpha'}\,([L_{12}]+[L_{23}]+[L_{13}])+
\frac{1}{2\alpha'}\,T,$$ where $1/2<\alpha'<\alpha$ is such that
$\nu(T,p_j)>(2-\alpha')/3$ for $j=1,2,3$. Then $\nu(R,q_j)>1/2$
and
$$\nu(R,p_j)>\frac{2(2\alpha'-1)}{6\alpha'}+\frac{2-\alpha'}{6\alpha'}
=\frac{1}{2}\;.$$ If $S=\{p_1,p_2,p_3,q_1,q_2\}$ then $m_1(S)=3$,
which contradicts Theorem \ref{T:234} $(ii)$.\end{pf}
\begin{Example}\label{E:magnet} Let $L_1,L_2$ be complex lines in
${\Bbb P}^n$ which intersect at the point $q$, and let
$T=\frac{1}{2}\,([L_1]+[L_2])$. Then $\nu(T,q)=1$ and
$E_{1/2}(T)\setminus L$ is uncountable for every complex line $L$.
So the assumption in Theorem \ref{T:magnet} on the existence of
two points where $T$ has large Lelong numbers, is
necessary.\end{Example}
\par Our last result shows that the complex lines from the
conclusion of Theorem \ref{T:234} are determined by 3 points where
$T$ has ``small" Lelong numbers.
\begin{Proposition}\label{P:24} Let $\alpha\geq1/2$ and assume
that $E_{1-\alpha}(T)$ contains the points $p_1,p_2,p_3$ which lie
on a complex line $L\subset{\Bbb P}^n$. Then
$|V_\alpha(T)\setminus L|\leq1$. Moreover, if $\alpha\geq2/3$ then
either $V_\alpha(T)\subseteq L$ or else $|V_\alpha(T)|\leq2$.
\end{Proposition}
\begin{pf} If $|V_\alpha(T)\setminus L|>1$, let $q_1,q_2\in
V_\alpha(T)\setminus L$, and let $\alpha'>\alpha$ be chosen so
that $\nu(T,q_j)>\alpha'$, $j=1,2$. Then
$$R=\frac{2\alpha'-1}{2\alpha'}\,[L]+\frac{1}{2\alpha'}\,T$$ has
Lelong numbers larger than $1/2$ at the points of the set
$S=\{p_1,p_2,p_3,q_1,q_2\}$. As $m_1(S)=3$, this contradicts
Theorem \ref{T:234} $(ii)$. \par We assume next that
$\alpha\geq2/3$ and there exists $q\in V_\alpha(T)\setminus L$. By
what we already proved, we have
$V_\alpha(T)\setminus\{q\}\subseteq L$. Theorem \ref{T:234} $(i)$
implies that $V_\alpha(T)\subseteq L'$, where $L'$ is a complex
line containing $q$. It follows that $|V_\alpha(T)|\leq2$.\end{pf}

\end{document}